\documentclass{llncs}
\usepackage{makeidx}

\usepackage{graphicx}
\usepackage{multirow}
\usepackage{color}
\usepackage{latexsym}
\usepackage{amssymb}
\usepackage{amsfonts}
\usepackage{amsmath}
\usepackage{indentfirst}
\usepackage{graphics}
\usepackage[customcolors]{hf-tikz}

\author{Elena Barcucci, Antonio Bernini, Renzo Pinzani}

\institute{
Dipartimento di Matematica e Informatica ``U. Dini'', Universit\`a degli
Studi di Firenze, Viale
 G.B. Morgagni 65, 50134 Firenze, Italy.\\  
\email{elena.barcucci@unifi.it, antonio.bernini@unifi.it, renzo.pinzani@unifi.it}
}

\title{Variable dimension non-overlapping 
	matrices}

\begin{document}
%
\pagestyle{headings} 
\maketitle

\begin{abstract}
	Since some years, non-overlapping sets of strings (also called cross-bifix-free sets) have had an increasing interest in the frame of the researches about Theory of Codes. Recently some non-overlapping sets of strings with variable length were introduced. Moreover, the notion of non-overlapping strings has been naturally extended to the two dimensional case leading to several definitions of non-overlapping sets of matrices (or pictures).
	
	Starting from these results, in this paper we introduce non-overlapping sets of binary matrices having variable dimension. Indeed, we use non-overlapping variable length strings as rows of the matrices and imposing the avoidance of two consecutive patterns of length $k$, we get the desired sets of non-ovelapping matrices with variable dimension.
	
\end{abstract}

\section{Introduction}
Given a finite alphabet $\Sigma$, a set $S\subset \Sigma^n$, containing strings of length $n$, is called \emph{cross-bifix-free} or \emph{non-overlapping} if any two strings $u,v\in S$ (the case $u=v$ is allowed) are such that any proper suffix of $u$ is different from any proper prefix of $v$ and vice versa (also the two strings $u,v$ are called non overlapping or cross-bifix-free). Several non-overlapping sets of strings have been defined in literature with different characteristics. For example, in \cite{BPP} the non-overlapping strings are binary strings and their definition involves Dyck paths, while in \cite{BBPPS} the alphabet is not necessarily a binary alphabet and the strings are defined via Motzkin paths. The set we can find in \cite{CKPW} is constituted by strings which avoid a particular consecutive pattern, while in \cite{Bl} the author investigates on the cardinality of certain sets depending on the length of the strings and the cardinality of the alphabet. In all the cited papers, the strings of the non-overlapping sets have the same length. Nevertheless it is possible to define non-overlapping strings with variable length \cite{Bi}. Herein, two different two different sets of non-overlapping strings are defined.
 
The notion of non-overlapping strings can be naturally 
extended to the two dimensional case by means of 
matrices (or pictures): any two of them $A,B$ (the case 
$A=B$ is allowed) do not overlap if it is not possible 
to move $A$ ($B$) over $B$ ($A$) in a way such that the 
corresponding entries match. Also in this case, it is possible to find different kinds of definitions of non-overlapping sets of matrices, as in \cite{BBBP1},\cite{BBBP2},\cite{BBBP3}, or several deep and interesting investigations about some their properties (\cite{AGM1,AGM2,AGM3,AGM4,AGM5} and all the references therein). In particular, in \cite{AGM1} the role of the frame of the matrices is deeply analysed, showing that it is possible to generate a corss-bifix-free set of matrices simply framing with a suitable frame any matrix of a given set. The approach used in\cite{BBBP1},\cite{BBBP2},\cite{BBBP3} is different: here the matrices are constructed by imposing some constraints on their rows which have to avoid some particular consecutive patterns or must have some fixed entries in particular positions.
The matrices of the sets defined in all the above cited papers have the same fixed dimension.

In this work, we move from variable length non-ovelapping sets of binary strings as defined in \cite{Bi} in order to define non-overlapping sets of binary matrices having variable dimension. Essentially, the strings are used to fill the rows of the matrices and, in a first approach, they have to avoid two particular consecutive patterns of length $k$. More precisely, fixed $m,n>0$, we define (non-overlapping) matrices having at most $m$ rows and at most $n$ columns. Actually, the number of columns must be at least $2k+3$. The cardinality of the sets of non-overlapping matrices defined in this first way is easily determined since it strictly depends on the cardinality of the used sets of strings \cite{Bi}.

Furthermore, we use a different set of variable length non-overlapping strings defined via Dyck paths. A set of non-overlapping matrices with variable dimension is defined and also in this case the cardinality is analysed.

\smallskip
Probably, the first natural application of non-overlapping sets of strings is in the theory of codes. We recall that a set of non-overlapping strings with the same length surely is a code (uniform code) with the further property that any two strings do not overlap. This fact has been used in the digital communication systems to establish and maintain a connection between a receiver and a transmitter. Secondly, the avoidance of overlaps in strings is a fundamental point in the context of string matching, compression and synchronization problems.
On the other side, the increasing interest for digital image processing is the reason why many works on the theory of two-dimensional languages have appeared, and where non-overlapping sets of matrices are framed. Maybe, a possible (future) application of this kind of sets is in the template matching which is a technique to discover if  small parts of an image match a template image.

\section{Preliminaries}

We first briefly recall the construction of the set of non-overlapping strings with variable length \cite{Bi} which is the starting point for the construction of our set of non-overlapping matrices with variable dimension.

\begin{definition}\label{strings}
Two strings $u$, $v$, possibly the same, are said \emph{non-overlapping} (or \emph{cross-bifix-free}) if any proper prefix of $u$ is different from any proper suffix of $v$, and vice versa.
\end{definition}

In the case $u=v$, the string  $u$ is said \emph{self non-overlapping} (or \emph{bifix-free}).

\begin{definition}\label{sets}
	A set of strings is said \emph{non-overlapping} (or \emph{cross-bifix-free}) if any two elements of the set are non-overlapping strings.
\end{definition}

\bigskip
Let $\Sigma=\{0,1\}$ be the alphabet we are going to use. We also fix $k\geq3$. With $0^k$ and $1^k$ we indicate the strings constituted by $k$ consecutive $0$'s and $k$ consecutive $1$'s, respectively.
With $Av(0^k,1^k)$ we denote the set of binary strings avoiding the forbidden patterns $0^k$ and $1^k$, and $|u|$ denotes the length of a binary string $u$.

We consider the set $V^{i,(k)}$ of strings with length $i\geq 2k+2$:

$$
V^{i,(k)}=\{1^k0u10^k|0u1\in Av(0^k,1^k),|u|\geq 0\}\ ,
$$

\noindent
 containing strings $v=v_1v_2\ldots v_i$  starting with $k$ consecutive $1$'s followed by a $0$, ending with $1$ followed by $k$ consecutive $0$'s, and avoiding $0^k$ and $1^k$ in their inner part (from the $k+1$-th letter to the $k+1$-th to last letter).
 
 Fixed $n\geq 2k+2$, the strings with length at most $n$ are collected in a set:
 
$$
\mathcal V_n^{(k)}=\bigcup_{i=2k+2}^{n}V^{i,(k)}\ .
$$
\noindent
Clearly, the strings in $\mathcal V_n^{(k)} $ have different lengths as soon as $n\gneqq 2k+2$. The following proposition holds:

\begin{proposition}
	The set $\mathcal V_n^{(k)}$ is a non-overlapping set of strings.
\end{proposition}

In Figure \ref{V_{13}^{(3)}} the set $\mathcal V_{13}^{(3)}$ containing the strings with length from $8$ to $13$ are showed. 	

\begin{figure}
\begin{center}
	{\footnotesize
		\begin{tabular}{|l|c||l|c|}
			\hline
			$n=8$	&11101000	& 	$n=12$		& 111010011000\\
			\cline{1-2}				
			$n=9$	&111011000	&	& 111011011000\\
			&111001000	&   & 111011001000\\
			\cline{1-2}				
			$n=10$ 	&1110101000 &   & 111010101000\\
			&1110011000 &   & 111010011000\\
			\cline{1-2}
			$n=11$ 	&11101101000&	&111001101000\\
			&11101001000&	&111001001000\\
			\cline{3-4}		
			&11101011000&$n=13$&1110011011000\\
			&11100101000& 	&1110010011000	\\
			\cline{1-2}
			&&&1110011001000\\
			&&&1110010101000\\
			&&&1110110011000\\
			&&&1110110101000\\
			&&&1110100101000\\
			&&&1110101101000\\
			&&&1110101011000\\
			&&&1110101001000\\
			\hline	
		\end{tabular}
	}
\end{center}
\caption{The set $\mathcal V_{13}^{(3)}$}
\label{V_{13}^{(3)}}
\end{figure}

An interesting property of the strings belonging to $\mathcal V_n^{(k)}$ is contained in the following corollary:

\begin{corollary}
	A string $u\in \mathcal V_n^{(k)}$  can not be a factor of a string $v\in \mathcal V_n^{(k)}$, for any two strings $u,v$.
\end{corollary}	 

Note that Definition \ref{strings} does not forbid the case excluded by the Corollary. For example $u=1110110000$ and $v=11011000$ are two non-overlapping strings, but $u=1v0$.

\section{Non-overlapping matrices}
 A rigorous definition of overlapping matrices can be 
 found in \cite{AGM1,AGM4} or in \cite{BBBP2}. Nevertheless, for our purpose, it 
 is enough to describe them as follows:
 given two matrices over an alphabet $\Sigma$, we start from a configuration where they are placed in a way such that their top left corners coincide. By moving the smaller one down, or right, or both, if the entries of the smaller matrix coincide with the one of the bigger matrix, we have two \emph{overlapping matrices}. If not, they are \emph{non-overlapping matrices}.
 
 In Figure \ref{overlapping matrices} some couples of matrices (with the same dimension) are presented which overlaps in three different ways. In the first case we say that they present a \emph{horizontal overlap}, in the second one they present a \emph{diagonal overlap}, finally in the third case they present a \emph{vertical overlap}.
 
  \begin{figure}
   	\begin{center}
   $$
 \begin{tikzpicture}[scale=0.30]
 \draw (0,0.1) rectangle (6,4.1);
 \draw (4,0) rectangle (10,4);
 \node at (0.5,3.5) {${0}$};
 \node at (1.5,3.5) {${1}$};
 \node at (2.5,3.5) {${0}$};
 \node at (3.5,3.5) {${1}$};
 \node at (4.5,3.5) {$\mathbf{0}$};
 \node at (5.5,3.5) {$\mathbf{1}$};
 
 \node at (0.5,2.5) {${0}$};
 \node at (1.5,2.5) {${1}$};
 \node at (2.5,2.5) {${1}$};
 \node at (3.5,2.5) {${1}$};
 \node at (4.5,2.5) {$\mathbf{0}$};
 \node at (5.5,2.5) {$\mathbf{1}$};
 
 \node at (0.5,1.5) {${1}$};
 \node at (1.5,1.5) {${1}$};
 \node at (2.5,1.5) {${1}$};
 \node at (3.5,1.5) {${1}$};
 \node at (4.5,1.5) {$\mathbf{0}$};
 \node at (5.5,1.5) {$\mathbf{0}$};
 
 \node at (0.5,0.5) {${1}$};
 \node at (1.5,0.5) {${0}$};
 \node at (2.5,0.5) {${0}$};
 \node at (3.5,0.5) {${0}$};
 \node at (4.5,0.5) {$\mathbf{1}$};
 \node at (5.5,0.5) {$\mathbf{1}$};
 

 \node at (6.5,3.5) {${1}$};
 \node at (7.5,3.5) {${1}$};
 \node at (8.5,3.5) {${1}$};
 \node at (9.5,3.5) {${0}$};

 \node at (6.5,2.5) {${1}$};
 \node at (7.5,2.5) {${1}$};
 \node at (8.5,2.5) {${1}$};
 \node at (9.5,2.5) {${0}$};
 
 \node at (6.5,1.5) {${1}$};
 \node at (7.5,1.5) {${0}$};
 \node at (8.5,1.5) {${0}$};
 \node at (9.5,1.5) {${0}$};

 \node at (6.5,0.5) {${0}$};
 \node at (7.5,0.5) {${1}$};
 \node at (8.5,0.5) {${0}$};
 \node at (9.5,0.5) {${1}$};
 \end{tikzpicture}
 $$
 
 $$
 \begin{tikzpicture}[scale=0.30]
 \draw (0,0) rectangle (6,4);
 \draw (3,-2) rectangle (9,2);
 \node at (0.5,3.5) {${0}$};
 \node at (1.5,3.5) {${1}$};
 \node at (2.5,3.5) {${0}$};
 \node at (3.5,3.5) {${1}$};
 \node at (4.5,3.5) {${0}$};
 \node at (5.5,3.5) {${1}$};
 
 \node at (0.5,2.5) {${0}$};
 \node at (1.5,2.5) {${1}$};
 \node at (2.5,2.5) {${1}$};
 \node at (3.5,2.5) {${1}$};
 \node at (4.5,2.5) {${0}$};
 \node at (5.5,2.5) {${1}$};
 
 \node at (0.5,1.5) {${1}$};
 \node at (1.5,1.5) {${1}$};
 \node at (2.5,1.5) {${1}$};
 \node at (3.5,1.5) {$\mathbf{1}$};
 \node at (4.5,1.5) {$\mathbf{0}$};
 \node at (5.5,1.5) {$\mathbf{0}$};
 
 \node at (0.5,0.5) {${1}$};
 \node at (1.5,0.5) {${0}$};
 \node at (2.5,0.5) {${0}$};
 \node at (3.5,0.5) {$\mathbf{0}$};
 \node at (4.5,0.5) {$\mathbf{1}$};
 \node at (5.5,0.5) {$\mathbf{1}$};
 
 \node at (6.5,0.5) {${1}$};
 \node at (7.5,0.5) {${1}$};
 \node at (8.5,0.5) {${0}$};
 \node at (6.5,1.5) {${1}$};
 \node at (7.5,1.5) {${1}$};
 \node at (8.5,1.5) {${1}$};
 
 \node at (3.5,-0.5) {${1}$};
 \node at (4.5,-0.5) {${1}$};
 \node at (5.5,-0.5) {${1}$};
 \node at (6.5,-0.5) {${0}$};
 \node at (7.5,-0.5) {${0}$};
 \node at (8.5,-0.5) {${0}$};
 
 \node at (3.5,-1.5) {${0}$};
 \node at (4.5,-1.5) {${1}$};
 \node at (5.5,-1.5) {${0}$};
 \node at (6.5,-1.5) {${1}$};
 \node at (7.5,-1.5) {${0}$};
 \node at (8.5,-1.5) {${1}$};
 \end{tikzpicture}
 $$
 
 $$
 \begin{tikzpicture}[scale=0.30]
 \draw (0,0) rectangle (6,4);
 \draw (-0.1,-2) rectangle (5.9,2);
 \node at (0.5,3.5) {${0}$};
 \node at (1.5,3.5) {${1}$};
 \node at (2.5,3.5) {${0}$};
 \node at (3.5,3.5) {${1}$};
 \node at (4.5,3.5) {${0}$};
 \node at (5.5,3.5) {${1}$};
 
 \node at (0.5,2.5) {${0}$};
 \node at (1.5,2.5) {${1}$};
 \node at (2.5,2.5) {${1}$};
 \node at (3.5,2.5) {${1}$};
 \node at (4.5,2.5) {${0}$};
 \node at (5.5,2.5) {${1}$};
 
 \node at (0.5,1.5) {$\mathbf{1}$};
 \node at (1.5,1.5) {$\mathbf{1}$};
 \node at (2.5,1.5) {$\mathbf{1}$};
 \node at (3.5,1.5) {$\mathbf{1}$};
 \node at (4.5,1.5) {$\mathbf{0}$};
 \node at (5.5,1.5) {$\mathbf{0}$};
 
 \node at (0.5,0.5) {$\mathbf{1}$};
 \node at (1.5,0.5) {$\mathbf{0}$};
 \node at (2.5,0.5) {$\mathbf{0}$};
 \node at (3.5,0.5) {$\mathbf{1}$};
 \node at (4.5,0.5) {$\mathbf{1}$};
 \node at (5.5,0.5) {$\mathbf{0}$};
 
 \node at (0.5,-0.5) {${1}$};
 \node at (1.5,-0.5) {${1}$};
 \node at (2.5,-0.5) {${1}$};
 \node at (3.5,-0.5) {${0}$};
 \node at (4.5,-0.5) {${0}$};
 \node at (5.5,-0.5) {${0}$};
 
 \node at (1.5,-1.5) {${1}$};
 \node at (0.5,-1.5) {${1}$};
 \node at (2.5,-1.5) {${0}$};
 \node at (3.5,-1.5) {${1}$};
 \node at (4.5,-1.5) {${0}$};
 \node at (5.5,-1.5) {${1}$};
 
 \end{tikzpicture}
 $$
\end{center}
\caption{Two horizontal (upper figure), diagonal (middle figure), and vertical (lower figure) overlapping matrices}
\label{overlapping matrices}
\end{figure}

As an example, the following to matrices are non-overlapping:

$$\begin{pmatrix}
1&1&0&0&1&0&0&1&0&0\\
0&1&1&0&0&1&0&1&1&0\\
0&1&0&0&1&1&0&0&1&0\\
1&1&1&1&1&0&0&0&0&0\\
\end{pmatrix}
\ \ \ \
\begin{pmatrix}
1&1&0&0&1&0&0&1&0&0\\
0&0&1&0&1&1&0&1&1&0\\
1&1&0&1&1&0&1&0&1&0\\
1&1&1&0&0&0&0&0&0&0\\
\end{pmatrix}\ .$$

\bigskip

We now define a set of variable dimension matrices, using strings of a same length $s$ of $V^{s,(k)}$ as rows of a matrix. In the following, the two matrices $C$ and $D$ of dimension $m_1\times s$ and $m_2\times t$, respectively, are constructed with the rows $C_i^{s,(k)}\in V^{s,(k)}$ and $D_j^{t,(k)}\in V^{t,(k)}$, with $i=1,2,\ldots,m_1$ and $j=1,2,\ldots, m_2$. 
$$
C=
\left(
\begin{matrix}
C_1^{s,(k)}\\
C_2^{s,(k)}\\
\vdots\\
\vdots\\
C_{m_1}^{s,(k)}
\end{matrix}
\right)
\hspace{1cm}
D=
\left(
\begin{matrix}
D_1^{t,(k)}\\
D_2^{t,(k)}\\
\vdots\\
D_{m_2}^{t,(k)}
\end{matrix}
\right)
$$

It is not difficult to show that $C$ and $D$ can present neither diagonal overlap nor horizontal overlap. If this were not the case, that is, if $C$ and $D$ had a horizontal or diagonal overlap, then there would be a prefix or suffix of some $D_j^{t,(k)}$ overlapping with some suffix or prefix of some $C_i^{t,(k)}$, against the hypothesis that $\mathcal V_n^{(k)}$ is a non-overlapping set.

Unfortunately, in the case $C$ and $D$ have the same number of columns ($s=t$), then the two matrices can present a vertical overlap, as showed in the following example:

$$C=
\left(
\begin{array}{cccccccccccc}
1&1&1&0&1&0&0&1&1&0&0&0\\
1&1&1&0&1&1&0&1&1&0&0&0\\
\textbf{1}&\textbf{1}&\textbf{1}&\textbf{0}&\textbf{1}&\textbf{1}&\textbf{0}&\textbf{0}&\textbf{1}&\textbf{0}&\textbf{0}&\textbf{0}\\
\textbf{1}&\textbf{1}&\textbf{1}&\textbf{0}&\textbf{1}&\textbf{0}&\textbf{1}&\textbf{0}&\textbf{1}&\textbf{0}&\textbf{0}&\textbf{0}\\
\textbf{1}&\textbf{1}&\textbf{1}&\textbf{0}&\textbf{1}&\textbf{0}&\textbf{0}&\textbf{1}&\textbf{1}&\textbf{0}&\textbf{0}&\textbf{0}\\
1&1&1&0&0&1&1&0&1&0&0&0\\
1&1&1&0&0&1&0&0&1&0&0&0\\
\end{array}
\right)
$$

\bigskip

\bigskip
$$D=
\left(
\begin{array}{cccccccccccc}
1&1&1&0&1&1&0&0&1&0&0&0\\
1&1&1&0&1&0&1&0&1&0&0&0\\
1&1&1&0&1&0&0&1&1&0&0&0\\
\end{array}
\right)
=\begin{pmatrix}
C_3^{12,(3)}\\
\\
C_4^{12,(3)}\\
\\
C_5^{12,(3)}\\
\end{pmatrix}
$$

It may also happens that a matrix $D$ has the first $h$ rows equal to the last $h$ rows of $C$, so that a vertical overlapping occurs again.
In order to avoid these kinds of overlaps we put a constraint in the first and in the last row of each matrix. In particular, all the matrices with the same number of columns must have the same first row and the same last row, different one from each other. Also, these two selected rows cannot appear as inner rows of any other matrix with that number of columns. In other words, we force:
\begin{itemize}
	\item the top and bottom row (T and B, respectively) of all the matrices with the same number of columns to be  the same for each matrix; 
	\item the rows $T$ and $B$ not to occur in any other line of the matrix.
\end{itemize} 

More precisely, the matrices $C$ with the same number $s$ of columns  must have the following structures:

$$
C=
\left(
\begin{matrix}
T^{s,(k)}\\
\\
C_2^{s,(k)}\\
\vdots\\
\vdots\\
C_{m_1-1}^{s,(k)}\\
\\
B^{s,(k)}
\end{matrix}
\right)
$$

\bigskip

\noindent
with $C_j^{s,(k)}\neq T^{s,(k)},B^{s,(k)}$,
for
$j= 2,3,\ldots, m_1-1$, and $C_j^{s,(k)},T^{s,(k)},B^{s,(k)}\in V^{s,(k)}$.
It is not difficult to convince ourselves that fixing the first and the last rows do not allow the occurrences of the vertical overlapping described above. Moreover, in order to construct a matrix with the strings of $V^{i,(k)}$, it must be $|V^{i,(k)}|\geq 2$. So, $i\geq 2k+3$. 

From the above discussion, it appears that there are no constraints in the choice of $T^{s,(k)}$ and $B^{s,(k)}$ among all the possible strings in $V^{s,(k)}$. For an example, we present the case  
$T^{s,(k)}=1^k0u10^k$ with $u=1010\ldots$ having suitable length, and 
$B^{s,(k)}=1^k0v10^k$ with $v=0101\ldots$ having suitable length, or, more precisely, $u=(10)^{\frac{s-2k-2}{2}}$ and $v=(01)^\frac{s-2k-2}{2}$ if $s$ is even and  $u=1(01)^{\frac{s-2k-3}{2}}$ and $v=0(10)^\frac{s-2k-3}{2}$, if $s$ is odd. In Figure \ref{esempio}, we show some non-overlapping matrices with $k=3$ and number of columns $n=9,10,11,12,13$. Note that in each group, the first (last) row of each matrix is the same, according to our particular choice. Moreover, the inner rows can be repeated within the same matrix.

\begin{figure}
	\begin{center}

\arraycolsep=0.78\arraycolsep
$n=9\quad\  
\begin{pmatrix}
1&1&1&0&1&1&0&0&0\\
1&1&1&0&0&1&0&0&0\\
\end{pmatrix}$

\medskip
$n=10\quad
\begin{pmatrix}
1&1&1&0&1&0&1&0&0&0\\
1&1&1&0&0&1&1&0&0&0\\
\end{pmatrix}$

\medskip
\arraycolsep=0.78\arraycolsep
$n=11
\quad
\left(
\begin{array}{ccccccccccc}
\textbf{1}&\textbf{1}&\textbf{1}&\textbf{0}&\textbf{1}&\textbf{0}&\textbf{1}&\textbf{1}&\textbf{0}&\textbf{0}&\textbf{0}\\
\textbf{1}&\textbf{1}&\textbf{1}&\textbf{0}&\textbf{0}&\textbf{1}&\textbf{0}&\textbf{1}&\textbf{0}&\textbf{0}&\textbf{0}\\
\end{array}
\right)
,
\left(
\begin{array}{ccccccccccc}
\textbf{1}&\textbf{1}&\textbf{1}&\textbf{0}&\textbf{1}&\textbf{0}&\textbf{1}&\textbf{1}&\textbf{0}&\textbf{0}&\textbf{0}\\
1&1&1&0&0&1&0&1&0&0&0\\
1&1&1&0&1&1&0&1&0&0&0\\
\textbf{1}&\textbf{1}&\textbf{1}&\textbf{0}&\textbf{0}&\textbf{1}&\textbf{0}&\textbf{1}&\textbf{0}&\textbf{0}&\textbf{0}\\
\end{array}
\right)
$

\medskip
$\qquad\qquad
\left(
\begin{array}{ccccccccccc}
\textbf{1}&\textbf{1}&\textbf{1}&\textbf{0}&\textbf{1}&\textbf{0}&\textbf{1}&\textbf{1}&\textbf{0}&\textbf{0}&\textbf{0}\\
1&1&1&0&1&0&1&1&0&0&0\\
\textbf{1}&\textbf{1}&\textbf{1}&\textbf{0}&\textbf{0}&\textbf{1}&\textbf{0}&\textbf{1}&\textbf{0}&\textbf{0}&\textbf{0}\\
\end{array}
\right)
,
\left(
\begin{array}{ccccccccccc}
\textbf{1}&\textbf{1}&\textbf{1}&\textbf{0}&\textbf{1}&\textbf{0}&\textbf{1}&\textbf{1}&\textbf{0}&\textbf{0}&\textbf{0}\\
1&1&1&0&0&1&0&1&0&0&0\\
1&1&1&0&1&1&0&1&0&0&0\\
1&1&1&0&1&1&0&1&0&0&0\\
\textbf{1}&\textbf{1}&\textbf{1}&\textbf{0}&\textbf{0}&\textbf{1}&\textbf{0}&\textbf{1}&\textbf{0}&\textbf{0}&\textbf{0}\\
\end{array}
\right),
\ \dots\dots$

\medskip

\medskip
\arraycolsep=0.78\arraycolsep
$n=12
\quad
\left(
\begin{array}{cccccccccccc}
\textbf{1}&\textbf{1}&\textbf{1}&\textbf{0}&\textbf{1}&\textbf{0}&\textbf{1}&\textbf{0}&\textbf{1}&\textbf{0}&\textbf{0}&\textbf{0}\\
1&1&1&0&0&1&0&1&1&0&0&0\\
1&1&1&0&0&1&0&0&1&0&0&0\\
\textbf{1}&\textbf{1}&\textbf{1}&\textbf{0}&\textbf{0}&\textbf{1}&\textbf{0}&\textbf{1}&\textbf{1}&\textbf{0}&\textbf{0}&\textbf{0}\\
\end{array}
\right)
,
\left(
\begin{array}{cccccccccccc}
\textbf{1}&\textbf{1}&\textbf{1}&\textbf{0}&\textbf{1}&\textbf{0}&\textbf{1}&\textbf{0}&\textbf{1}&\textbf{0}&\textbf{0}&\textbf{0}\\
1&1&1&0&0&1&0&0&1&0&0&0\\
1&1&1&0&1&1&0&0&1&0&0&0\\
\textbf{1}&\textbf{1}&\textbf{1}&\textbf{0}&\textbf{0}&\textbf{1}&\textbf{0}&\textbf{1}&\textbf{1}&\textbf{0}&\textbf{0}&\textbf{0}\\
\end{array}
\right)
$

\medskip
$\qquad\qquad
\left(
\begin{array}{cccccccccccc}
\textbf{1}&\textbf{1}&\textbf{1}&\textbf{0}&\textbf{1}&\textbf{0}&\textbf{1}&\textbf{0}&\textbf{1}&\textbf{0}&\textbf{0}&\textbf{0}\\
1&1&1&0&1&0&0&1&1&0&0&0\\
\textbf{1}&\textbf{1}&\textbf{1}&\textbf{0}&\textbf{0}&\textbf{1}&\textbf{0}&\textbf{1}&\textbf{1}&\textbf{0}&\textbf{0}&\textbf{0}\\
\end{array}
\right)
,
\left(
\begin{array}{cccccccccccc}
\textbf{1}&\textbf{1}&\textbf{1}&\textbf{0}&\textbf{1}&\textbf{0}&\textbf{1}&\textbf{0}&\textbf{1}&\textbf{0}&\textbf{0}&\textbf{0}\\
1&1&1&0&0&0&1&0&1&0&0&0\\
1&1&1&0&1&1&0&0&1&0&0&0\\
1&1&1&0&1&1&0&0&1&0&0&0\\
\textbf{1}&\textbf{1}&\textbf{1}&\textbf{0}&\textbf{0}&\textbf{1}&\textbf{0}&\textbf{1}&\textbf{1}&\textbf{0}&\textbf{0}&\textbf{0}\\
\end{array}
\right),
\ \dots\dots$

\medskip

$n=13
\quad
\left(
\begin{array}{ccccccccccccc}
\textbf{1}&\textbf{1}&\textbf{1}&\textbf{0}&\textbf{1}&\textbf{0}&\textbf{1}&\textbf{0}&\textbf{1}&\textbf{1}&\textbf{0}&\textbf{0}&\textbf{0}\\
1&1&1&0&0&1&0&0&1&1&0&0&0\\
1&1&1&0&0&1&0&1&1&1&0&0&0\\
1&1&1&0&0&1&0&0&1&1&0&0&0\\
\textbf{1}&\textbf{1}&\textbf{1}&\textbf{0}&\textbf{0}&\textbf{1}&\textbf{0}&\textbf{1}&\textbf{0}&\textbf{1}&\textbf{0}&\textbf{0}&\textbf{0}\\
\end{array}
\right)
,
\left(
\begin{array}{ccccccccccccc}
\textbf{1}&\textbf{1}&\textbf{1}&\textbf{0}&\textbf{1}&\textbf{0}&\textbf{1}&\textbf{0}&\textbf{1}&\textbf{1}&\textbf{0}&\textbf{0}&\textbf{0}\\
1&1&1&0&0&1&1&0&0&1&0&0&0\\
1&1&1&0&0&1&1&0&0&1&0&0&0\\
\textbf{1}&\textbf{1}&\textbf{1}&\textbf{0}&\textbf{0}&\textbf{1}&\textbf{0}&\textbf{1}&\textbf{0}&\textbf{1}&\textbf{0}&\textbf{0}&\textbf{0}\\
\end{array}
\right)
$

\medskip
$\qquad\qquad
\left(
\begin{array}{ccccccccccccc}
\textbf{1}&\textbf{1}&\textbf{1}&\textbf{0}&\textbf{1}&\textbf{0}&\textbf{1}&\textbf{0}&\textbf{1}&\textbf{1}&\textbf{0}&\textbf{0}&\textbf{0}\\
1&1&1&0&1&0&1&0&1&1&0&0&0\\
1&1&1&0&1&1&0&1&0&1&0&0&0\\
1&1&1&0&1&1&0&1&0&1&0&0&0\\
\textbf{1}&\textbf{1}&\textbf{1}&\textbf{0}&\textbf{0}&\textbf{1}&\textbf{0}&\textbf{1}&\textbf{0}&\textbf{1}&\textbf{0}&\textbf{0}&\textbf{0}\\
\end{array}
\right)
,
\left(
\begin{array}{ccccccccccccc}
\textbf{1}&\textbf{1}&\textbf{1}&\textbf{0}&\textbf{1}&\textbf{0}&\textbf{1}&\textbf{0}&\textbf{1}&\textbf{1}&\textbf{0}&\textbf{0}&\textbf{0}\\
1&1&1&0&1&0&0&1&0&1&0&0&0\\
1&1&1&0&1&1&0&1&0&1&0&0&0\\
1&1&1&0&1&0&0&1&0&1&0&0&0\\
1&1&1&0&1&0&0&1&0&1&0&0&0\\
\textbf{1}&\textbf{1}&\textbf{1}&\textbf{0}&\textbf{0}&\textbf{1}&\textbf{0}&\textbf{1}&\textbf{0}&\textbf{1}&\textbf{0}&\textbf{0}&\textbf{0}\\
\end{array}
\right),
\ \dots\dots
$

\end{center}
\caption{Some non-overlapping matrices with $k=3$ and variable dimension}
\label{esempio}
\end{figure}

Collecting our results, we define the set $\mathcal V_{m,n}^{(k)}$ of non-overlapping matrices having at most $n$ columns and at most $m$ rows:

\begin{definition}
	Let $$
	\mathcal V_{m,n}^{(k)}=
	\bigcup_
	{2\le h\le m \atop 2k+3\le s\leq n}
	M_{h,s}^{(k)}
	$$
	where
\tiny
	$$
	M_{h,s}^{(k)}=
	\left\{
	\begin{pmatrix}
	T^{s,(k)}\\
	A_2^{s,(k)}\\
	\vdots\\
	A_{h-1}^{s,(k)}\\
	B^{s,(k)}		
	\end{pmatrix}
	\middle|\ 
	A_i^{s,(k)},T^{s,(k)},B^{s,(k)}\in V^{s,(k)}
	\ \mbox{and}\ 
	A_i^{s,(k)}\neq T^{s,(k)},B^{s,(k)}\ \mbox{for}\ 
	i=2,3,\ldots, h-1
	\right\}
	$$
\end{definition}

\bigskip
For the above arguments, we have the following proposition:
\begin{proposition}
	The set $\mathcal V_{m,n}^{(k)}$
	is non-overlapping.	
\end{proposition}

\section{Cardinality}
The number of matrices in $\mathcal V_{m,n}^{(k)}$ strictly depends on the number of rows different from the first and the last one.
Each inner row is of the form $1^k0u10^k$, where $0u1$ is a binary string avoiding $0^k$ and $1^k$. Denoting by $R_{\ell}(0^k,1^k)$ the set of binary strings starting with $0$, ending with $1$, and avoiding $0^k$ or $1^k$, and with $r_{\ell}^{(k)}$ its cardinality, from the construction of the set $\mathcal V_{m,n}^{(k)}$, we deduce that its cardinality is:

{\footnotesize
\begin{equation}\label{eqn}
|\mathcal V_{m,n}^{(k)}|=\sum_{h=2}^{m}\sum_{s=2k+3}^{n}|M_{h,s}|
=\sum_{h=2}^{m}\sum_{s=2k+3}^{n}\left(r_{s-2k}^{(k)}-2\right)^{h-2}
=\sum_{h=2}^{m}\sum_{s=3}^{n-2k}\left(r_{s}^{(k)}-2\right)^{h-2}
\ .
\end{equation}
}

It is known \cite{BBBP2,Be} that 

\begin{equation*}\label{rlfib}
r_{\ell}^{(k)}=\left \{
\begin{array}{ll}
1 & \mbox{if} \ {\ell} = 0\\
\\
\displaystyle\frac{f^{(k-1)}_{{\ell}-1} + d^{(k)}_{\ell}}{2} & \mbox{if} \ {\ell} \geq 1\ \ ,
\end{array}
\right.
\end{equation*}

where $f_{\ell}^{(k)}$ are the $k$-generalized Fibonacci numbers

\begin{equation*}\label{kbonaccistring}
f_{\ell}^{(k)}=\left \{
\begin{array}{ll}
2^\ell & \mbox{if} \ 0 \leq {\ell} < k-1\\
\\
\displaystyle \sum_{i=1}^k f_{{\ell}-i}^{(k)} & \mbox{if} \ {\ell} \geq k\ ,
\end{array}\right.
\end{equation*}

and

\begin{equation*}
d^{(k)}_{\ell}=
\left\{
\begin{array}{rl}
1& \mbox{if}\ ({\ell}\ \mbox{mod}\ k)=0\\
\\
-1& \mbox{if}\ ({\ell}\ \mbox{mod}\ k)=1\\
\\
0& \mbox{if}\ ({\ell}\ \mbox{mod}\ k)\geq 2\ .\\
\end{array}
\right.
\end{equation*}

The first terms of the last sequence are: 
$$\left\{d_{\ell}^{(k)}\right\}_{\ell \geq 0}=\Big\{1,-1,
\underbrace{0,0,\ldots,0}_{k-2},1,-1,
\underbrace{0,0,\ldots,0}_{k-2},1,-1,0,\ldots\Big\}\ ,
$$
so that $-1\leq d_{\ell}^{(k)}\leq 1$ and, from (\ref{eqn}), we deduce

$$
\sum_{h=2}^{m}\sum_{s=3}^{n-2k}
\left(\frac{f_{s-1}^{(k)}-5}{2}\right)^{(h-2)}
\leq
|\mathcal V_m^{n,(k)}|
\leq
\sum_{h=2}^{m}\sum_{s=3}^{n-2k}
\left(\frac{f_{s-1}^{(k)}-3}{2}\right)^{(h-2)}\ \ ,
$$
showing that the number of non-overlapping matrices we defined grows fast.

\bigskip
We conclude this section with a short discussion about the so called non-expandability. First of all we recall its definition in the framework of strings:

\begin{definition}
Let $X$ be a non-overlapping set of strings and $X_n$ be the subset of $X$ containing the strings having length at most $n$. For any fixed $n$, we say that $X_n$ is non-expandable if and only if for each self non-overlapping string $\varphi$ having length at most $n$, with $\varphi \notin X_n$, we have that $X_n\cup \{\varphi\}$ is not a non-overlapping set.
\end{definition}
The set $\mathcal V_n^{(k)}$ is not non-expandable. It is not difficult to check that for each string
 $w=1^{\left\lceil\frac{\ell}{2}\right\rceil}
0^{\left\lfloor\frac{\ell}{2}\right\rfloor}$, $\ell=2k,2k+1,\ldots,n$ the set $\mathcal V_n^{(k)} \cup \{w\}$ is still non-overlapping.

\bigskip
For what non-overlapping matrices is concerned, we have the following definition:
\begin{definition}
	Let $Y$ be a non-overlapping set of matrices and $Y_{m,n}$ be the subset of $Y$ containing the matrices having at most $m$ rows and $n$ columns. For any fixed $m,n$, we say that $Y_{m,n}$ is \emph{non-expandable} if and only if for each self non-overlapping matrix $U$ having at most $m$ rows and $n$ columns, with $U \notin Y_{m,n}$, we have that $Y_{m,n} \cup \{U\}$ is not a non-overlapping set.
\end{definition}
Also in the case of our set we have that $\mathcal V_{m,n}^{(k)}$ is not non-expandable. Indeed, for each matrix
$W=
\begin{pmatrix}
T^{s,(k)}\\
1^{\left\lceil\frac{s}{2}\right\rceil}
0^{\left\lfloor\frac{s}{2}\right\rfloor}
\end{pmatrix}
$,
with $2k+3\le s\le n$, the set $\mathcal V_{m,n}^{(k)} \cup \{W\}$ is still non-overlapping. However, in a possible application we believe that non-expandability is not decisive for a set of non-overlapping matrices. For example, in the case of template matching we believe that a large amount of template images could be more useful than having a template images set which is non-expandable. Furthermore, if a set of non-overlapping matrices is non-expandable, another problem immediately arises: which is the maximal set? Often, the answer to this matter is a hard challenge and seems to be so, even in the case of non-overlapping matrices.

\section{Further developments and conclusion}

On a closer inspection, the definition of the set $\mathcal V_{m,n}^{(k)}$ does not depend on the particular set $\mathcal V_n^{(k)}$ of non-overlapping strings we used. Actually, the crucial point is the constraint on the first and last rows which must be the same for all the matrices with the same number of columns. Therefore, once a particular set of variable length non-overlapping strings is defined, it is possible to generate a set of non-overlapping matrices with variable dimension with the same technique.

In \cite{Bi}, a different set of variable dimension strings constructed via Dyck paths is proposed. It is the set  $\mathcal D_n$:
$$\mathcal D_n=\bigcup_{i\geq 0} ^{
	\left\lfloor		
	\frac{n-2}{2}	
	\right\rfloor
}\left\{1\omega0: \omega\in D_{2i}\right\}
$$
where $D_{2i}$ is the set of Dyck paths of length $2i$.

Given $u,v\in \mathcal D_n$, this time $u$ can be an inner factor of $v$, as showed by the strings $u = 11101000$ and $v=11100\textbf{11101000}100=11100u100$.
Note that the definition of non-overlapping strings does not involve inner factors. So, if we admit this kind of overlapping, we have  \cite{Bi}

\begin{proposition}
$\mathcal D_n$ is a variable length non-overlapping set of
strings.
\end{proposition}

Starting from here and following the same steps of Section 3  we define the set $\mathcal D_{m,n}$ of non overlapping matrices with variable dimension having at most $m$ rows and $n$ columns:

\begin{definition}
	Let $$
	\mathcal D_{m,n}=
	\bigcup_
	{2\le h\le m \atop 2\le s\leq \left\lfloor		
		\frac{n}{2}	
		\right\rfloor}
	M_{h,2s}
	$$
	where
	$$
	M_{h,2s}=
	\left\{
	\begin{pmatrix}
	T^{2s}\\
	A_2^{2s}\\
	\vdots\\
	A_{h-1}^{2s}\\
	B^{2s}		
	\end{pmatrix}
	\right\}
	$$
	with
	\begin{itemize}
		\item $T^{2s}=1u0,\ B^{2s}=1v0,\ u,v\in D_{2s-2}$;
		\item $A_{i}^{2s}=1\omega 0,\ \omega\in D_{2s-2}$; 
		\item $A_{i}^{2s}\neq T^{2s},\ B^{2s}$;	\end{itemize}
\end{definition}

\medskip
The following proposition holds.
\begin{proposition}
	The set $\mathcal D_{m,n}$
	is non-overlapping.	
\end{proposition}

\medskip
We note that each matrix of the above set contains the same number of $1$'s and $0$'s. Moreover, its cardinality can be easily deduced from the fact that, denoting with $C_n$ the $n$-th Catalan number counting the Dyck paths of semilength $n$, the $i$-th inner rows of a matrix with $2s$ columns is 
enumerated by $C_{s-1}$.
Then:
$$
	|\mathcal D_{m,n}|=
	\sum_{h=2}^{m}
	\sum_{s=2}^{\left\lfloor\frac{n}{2}	
		\right\rfloor}\left(C_{s-1}-2\right)^{h}\ .
$$

\medskip
We conclude with a general consideration about the construction of our sets of non-overlapping matrices with variable dimension. We followed a different approach with respect to the one used in \cite{AGM1},\cite{AGM4} or in \cite{BBBP1},\cite{BBBP2},\cite{BBBP3}, as mentioned in Introduction.
Here, we simply use existing sets of particular strings which are non-overlapping and having variable lengths, and then we transfer these two characteristics to a new set of bidimensional objects (matrices), generalizing two typical concepts of linear structures (strings). In this way, the construction of the matrices is easy and straightforward, moreover it allows an easy enumeration. 

\medskip

As a further development, it could be interesting to investigate about the possibility to list the matrices in a Gray code sense, following the studies started in \cite{BBP},\cite{BBPV1},\cite{BBPV2},\cite{BBPSV1},\cite{BBPSV2}, where different Gray codes are defined for several set of strings. Clearly, the variable dimension of the matrices has to be considered in a possible definition of a suitable Gray code.



\newpage
\begin{thebibliography}{99}
\bibitem{AGM1} Anselmo, M., Giammarresi, D., Madonia, M.:
Sets of pictures avoiding overlaps.
Internat. J. Found. Comput. Sci. {\bf 30}(6 \& 7), 875--898 (2019).

\bibitem{AGM2} Anselmo, M., Giammarresi, D., Madonia, M.:
Deterministic and unambiguous fam-ilies  within  recognizable  two-dimensional  languages.
Fund. Inform. {\bf 98}(2–3), 143–-166  (2010).

\bibitem{AGM3} Anselmo, M., Giammarresi, D., Madonia, M.:
A  computational  model  for  tiling recognizable two-dimensional languages.
Theor. Comput. Sci. {\bf 410}(37), 3520-–3529 (2009).

\bibitem{AGM4} Anselmo, M., Giammarresi, D., Madonia, M.:
Non-expandable non-overlapping sets of pictures.
Theor. Comput. Sci. {\bf 657}, 127-–136 (2017).

\bibitem{AGM5} Anselmo, M., Giammarresi, D., Madonia, M.: Picture codes and deciphering delay.
Inf. Comput. {\bf 253} 358--370 (2017).

\bibitem{BBBP1} Barcucci, E., Bernini, A., Bilotta, S., Pinzani, R.: Cross-bifix-free sets in two dimensions.
Theoret. Comput. Sci. {\bf 664}, 29--38 (2015).	
	
\bibitem{BBBP2} Barcucci, E., Bernini, A., Bilotta, S., Pinzani, R.: Non-overlapping matrices.
Theoret. Comput. Sci. {\bf 658}, 36--45 (2017).

\bibitem{BBBP3} Barcucci, E., Bernini, A., Bilotta, S., Pinzani, R.: A 2D non-overlapping code over a q-ary alphabet.
Cryptogr. Commun. {\bf 10}(4), 667--683 (2018).
	
\bibitem{BBP} Barcucci, E., Bernini, A., Pinzani, R.:
A Gray code for a regular language.
In. Ferrari, L., Vamvakari, M.
GASCom 2018, CEUR Workshop Proceedings,
vol. 2113, pp. 87--93 (2018).

\bibitem{BBPPS} Barcucci, E., Bilotta, S., Pergola, E., Pinzani, R., Succi, J.:
Cross-bifix-free sets generation via Motzkin paths. RAIRO-Theor. Inf. Appl. {\bf 50}(1), 81--91 (2016).

\bibitem{Be} Bernini, A.:
Restricted binary strings and generalized Fibonacci numbers.
In: Dennunzio, A., Formenti, E., Manzoni, L., Porreca, A. E.
AUTOMATA 2017, LNCS, vol. 10248, pp. 32--43.     Springer International Publishing, Milan (2017).   

\bibitem{BBPV1}
Bernini, A., Bilotta, S., Pinzani, R., Vajnovszki, V.:
A Gray code for cross-bifix-free sets.
Math. Structures Comput. Sci. {\bf 27}(2), 184--196 (2017).

\bibitem{BBPV2}
Bernini, A., Bilotta, S., Pinzani, R., Vajnovszki, V.:
A trace partitioned Gray code for q-ary generalized Fibonacci strings.
J. Discrete Math. Sci. Cryptogr. {\bf 18}(6), 751--761 (2015).

\bibitem{BBPSV1}
Bernini, A., Bilotta, S., Pinzani, R., Sabri, A., Vajnovszki, V.:
Gray code orders for q-ary words avoiding a given factor.
Acta Inform. {\bf 52}(7-8), 573--592 (2015).

\bibitem{BBPSV2}
Bernini, A., Bilotta, S., Pinzani, R., Sabri, A., Vajnovszki, V.:
Prefix partitioned gray codes for particular cross-bifix-free sets.
Cryptogr. Commun. {\bf 6}(4) 359--369 (2014).

\bibitem{Bi} Bilotta, S.:
Variable-length non ovelapping codes.
IEEE Trans. Inform. Theory {\bf 63}(10),  6530--6537 (2017).

\bibitem{BPP} Bilotta, S., Pergola, E., Pinzani, R.:
A new approach to cross-bifix-free sets.
IEEE Trans. Inform. Theory {\bfseries 58}(10),  4058--4063 (2012).

\bibitem{Bl} Blackburn, S. R.:
Non-overlapping codes.
IEEE Trans. Inform. Theory {\bf 61}(9),  4890--4894 (2015).

\bibitem{CKPW} Chee, Y. M., Kiah, H. M., Purkayastha, P., Wang, C.:
Cross-bifix-free codes within a constant factor of optimality.
IEEE Trans. Inform. Theory. {\bf 59}(7) 4668--4674 (2013).

%
%
%
%
%
%
%
%
%
%
%
%
%
%
%
%
%
%
%
%
%
%
%


\end{thebibliography}
\end{document}